\documentclass{amsart} 

\usepackage{epsfig}
\usepackage{amsmath}

\newtheorem{theorem}{Theorem}[section]
\newtheorem{proposition}[theorem]{Proposition}
\newtheorem{definition}[theorem]{Definition}
\newtheorem{lemma}[theorem]{Lemma}
\newtheorem{corollary}[theorem]{Corollary}
\newtheorem{remark}[theorem]{Remark}
\newtheorem{example}[theorem]{Example}
\newtheorem{conjecture}[theorem]{Conjecture}
\newtheorem{problem}[theorem]{Problem}
\newtheorem{question}[theorem]{Question}

\newcommand{\id}{\mathrm{id}}
\newcommand{\inv}{\mathrm{inv}}
\newcommand{\poi}{\mathrm{Poin}}
\newcommand{\fix}{\mathrm{Fix}(\theta)}
\newcommand{\twistid}{\mathfrak{I}(\id)}
\newcommand{\C}{\mathbb{C}}
\newcommand{\Z}{\mathbb{Z}}
\newcommand{\N}{\mathbb{N}}
\newcommand{\Br}{\mathrm{Br}}
\newcommand{\twist}{\mathfrak{I}(\theta)}
\newcommand{\s}{\underline{s}}

\newcommand{\SSS}{\underline{S}}
\newcommand{\iot}{\iota(\theta)}
\begin{document}

\title{Twisted identities in Coxeter groups}  
\author{Axel Hultman}
\address{Department of Mathematics, KTH, SE-100 44 Stockholm, Sweden} 
\email{axel@math.kth.se} 

\begin{abstract}
Given a Coxeter system $(W,S)$ equipped with an involutive
automorphism $\theta$, the set of {\em twisted identities} is 
\[
\iot = \{\theta(w^{-1})w\mid w\in W\}.
\]
We point out how $\iot$ shows up in several contexts and prove
that if there is no $s\in S$ such that $s\theta(s)$ is of odd order
greater than $1$, then the Bruhat order on $\iot$ is a graded poset with rank
function $\rho$ given by halving the Coxeter length. Under the same
condition, it is shown that the order complexes of the open intervals
either are PL spheres or $\Z$-acyclic. In the general case,
contractibility is shown for certain classes of
intervals. Furthermore, we demonstrate that sometimes these posets are
not graded.

For the Poincar\'e series of $\iot$, i.e.\ its generating
function with respect to $\rho$, a factorisation phenomenon is discussed. 
\end{abstract} 

\maketitle
\section{Introduction}
Let $(W,S)$ be a Coxeter system with an involutive automorphism
$\theta$. A {\em twisted identity} is an element of the form
$\theta(w^{-1})w$ for $w\in W$. In other words, the set $\iot$
of twisted 
identities is the orbit of the identity element under the twisted
conjugation action of $W$ on itself. Our terminology appeared in
\cite{hultman1} and stems from the fact 
that when $\theta = \id$, the only twisted identity is the identity
element $e$. 

As will be demonstrated in Section
\ref{se:motivation}, the study of $\iot$ is motivated by its
appeareance in a variety of situations. For 
example, certain orbit decompositions of symmetric varieties have a
close relationship with the subposet $\Br(\iot)$ of the Bruhat order
on $W$ which is induced by $\iot$ \cite{RS1, RS2}. Also, the Bruhat
order on $W$ itself appears as a special case of $\Br(\iot)$; see
Example \ref{ex:cox}. 

This article is chiefly devoted to $\Br(\iot)$. The
questions that we strive to answer emerge from a context which we now
briefly describe.

A {\em twisted involution} in $W$ is an element which is sent to its
inverse by $\theta$. Denote by $\twist$ the set of twisted
involutions. Clearly, 
\[
\iot \subseteq \twist \subseteq W.
\]
Let $\Br(X)$ denote the subposet of the Bruhat order on $W$ induced by
$X\subseteq W$. It is a fact that $\Br(W)$ is graded with
rank function given by the Coxeter length $\ell$. Furthermore, a
fundamental result due to Bj\"orner and
Wachs \cite{BW} asserts that (the order complex of) every (open)
interval in $\Br(W)$ is a PL sphere. Recent results on 
$\Br(\twist)$ produce a similar picture. It is graded with rank
function $\rho = (\ell+\ell^\theta)/2$, where
$\ell^\theta$ is the {\em twisted absolute length} \cite{hultman1}. Moreover,
every interval in $\Br(\twist)$ is a PL sphere \cite{hultman2}.

In the spirit of the above description, we pose the
following problems:
\begin{itemize}
\item[1.] For which $(W,S)$ and $\theta$ is $\Br(\iot)$ a graded poset?
\item[2.] Describe the topology of the intervals of $\Br(\iot)$. 
\end{itemize}
We do not know the complete solution to either of the problems. Our
main results on $\Br(\iot)$ are these partial answers:
 
\begin{itemize}
\item In Theorem \ref{th:graded} it is shown that if $s\theta(s)$
  never is of odd order for $s\in S$ unless 
  $s$ is a fixed point of $\theta$, then $\Br(\iot)$ is graded with
  rank function $\rho$. For
  example, this is always the case if $W$ is of odd rank and its
  Coxeter graph is a tree. By way of contrast, there exist $(W,S)$ and
  $\theta$ such that $\Br(\iot)$ is not graded; see Example \ref{ex:ungraded}.
\item Under the same conditions on $s\theta(s)$ as above, every
  interval in $\Br(\iot)$ is either $\Z$-acyclic (i.e.\ has trivial
  reduced integral homology) or a PL sphere. The
  latter case occurs precisely when the interval coincides with an
  interval in $\Br(\twist)$. This is Theorem \ref{th:topology}. Dropping these
  conditions, the homotopy type of an interval can be
  computed in certain special cases (Theorems \ref{th:collapsible} and
  \ref{th:recursive}).
\end{itemize}

In addition to the above results, we also investigate what we call the
{\em Poincar\'e series} of $\iot$. This is the rank generating function
of $\Br(\iot)$ whenever it is graded. Specifically, we provide a simple
necessary condition for an intriguing factorisation phenomenon to
occur and demonstrate that this condition also is sufficient in the
context of finite Coxeter groups. The general case is left open.

The structure of the remainder of the paper is as follows. In the next
section we gather preliminaries on poset topology and Coxeter groups, including
some new material on twisted involutions, that we need in
the sequel. In Section \ref{se:motivation} we then give an account of
various natural contexts where $\iot$ and $\Br(\iot)$ appear. After that, we turn to the study of $\Br(\iot)$ in Section 
\ref{se:bruhat}. Section \ref{se:poincare} is concerned with
Poincar\'e series considerations. Finally, we mention several open
problems and give further comments in Section \ref{se:problems}.

\section{Preliminaries} \label{se:preliminaries}

\subsection{Posets and combinatorial topology}
Say that a poset is {\em locally finite} if every interval is
  finite. A locally finite poset equipped with a minimum element is
  {\em graded} if in any given interval all maximal chains have the same
  cardinality. 

With any finite poset $P$, we associate the {\em order complex}
$\Delta(P)$. It is the (abstract) simplicial complex whose simplices
are the chains in $P$. Whenever we make topological statements about
$P$, we have the corresponding properties of $\Delta(P)$ 
in mind. We make no notational distinction between
an abstract simplicial complex and its geometric realisation.

Given a regular CW complex $\Delta$, we define the {\em face poset}
$P(\Delta)$ to be the set of cells in $\Delta$ ordered by
inclusion. By convention, we include the empty cell as the minimum
element in $P(\Delta)$.

A simplicial complex is a {\em PL (piecewise linear) ball} if it has
a common subdivision with a simplex. Similarly, it is
a {\em PL sphere} if it has a common subdivision with the boundary of
a simplex. In particular, a PL sphere is of course homeomorphic to a
sphere (and similarly for balls).

We now review elements of Forman's discrete Morse theory
\cite{forman}. Its formulation in terms of matchings, to which we
adhere, is due to Chari \cite{chari}. For unexplained terminology from
combinatorial topology as well as further background, we refer to
\cite{bjorner}. 

Let $\Delta$ be a finite regular CW complex. A {\em matching} on the face
poset $P(\Delta)$ is an involution $M:Q \to Q$, for some $Q\subseteq
P(\Delta)$ such that for all $q\in Q$, either $M(q)\prec q$ or
$M(q)\succ q$, where $\prec$ is the covering relation in
$P(\Delta)$. In other words, $M$ is nothing but a graph-theoretic
matching of the Hasse diagram of $P(\Delta)$. The unmatched elements,
i.e.\ the members of $P(\Delta)\setminus Q$, are called the {\em
  critical cells}.

The matching $M$ is {\em acyclic} if whenever we have a sequence in
$Q$ of the form 
\[
q_0 \prec M(q_0) \succ q_1 \prec M(q_1) \succ \cdots \prec M(q_{t-1})
\succ q_t
\]
with $q_1 \neq q_0$, it holds that $q_t \neq q_0$. A nice
way to interpret this condition is as follows. If, in the Hasse
diagram of $\Delta(P)$, we direct the matching edges upwards and the
others downwards, then this directed graph is acyclic iff $M$ is an
acyclic matching.

Without loss of generality, we will always assume that the empty cell is not
critical, i.e.\ that $Q$ includes the minimum
element\footnote{Otherwise, we may just extend $M$ to include
  it. There is always a matching partner for the empty cell available, 
because no acyclic matching can match all $0$-dimensional cells with
$1$-dimensional ones.}.

Our acyclic matching $M$ determines a way to collapse $\Delta$ onto a
(possibly non-regular) {\em critical complex} $\Delta^M$
consisting of the critical cells 
together with the vertex which was matched with the empty cell. In the
process, incidences among the cells may change. In some
situations, however, this is not a problem.

\begin{theorem}[Forman \cite{forman}] \label{th:forman}
If $M$ is an acyclic matching on $P(\Delta)$ which is complete (i.e.\
has no critical cells), then $\Delta$ is collapsible, in particular
contractible. 
\end{theorem}

Theorem \ref{th:forman} has the following consequence, which happens
to be well-suited for some of our applications.

\begin{corollary} \label{co:filter}
Suppose $M$ is an acyclic matching on $P(\Delta)$ with set of
critical cells $C$. If $C$ is an order filter (i.e.\ $c\in C$, $b\geq
c \Rightarrow b\in C$), then $\Delta$ is homotopy equivalent to the
complex $\Delta^M$ obtained from $\Delta$ by collapsing the complex of
non-critical cells to a point.  
\begin{proof}
A standard result in topology asserts that for a CW
complex $A$ with contractible subcomplex $A_0$, the quotient map $A\to
A/A_0$ is a homotopy equivalence. Since $C$ is an order filter,
$\Delta \setminus C$ is a subcomplex of $\Delta$. This subcomplex is
contractible by Theorem \ref{th:forman}.
\end{proof}
\end{corollary}

\subsection{Coxeter groups}
The reader is assumed to be familiar with the basics of Coxeter group
theory. Here, bits and pieces are reviewed in order to agree on
notation. We refer to the textbooks \cite{BB} and \cite{humphreys} for
further details.

Henceforth, let $(W,S)$ be a Coxeter system with $|S|<\infty$. The
Coxeter length function is $\ell:W\to \N$. If $w = s_1\cdots s_k \in W$
and $\ell(w) = k$, the word $s_1\cdots s_k$ is called a {\em reduced
  expression} for $w$. Here, and in what follows, symbols of the form
$s_i$ are always elements in $S$. We make no notational distinction
between a word in the free monoid over $S$ and the element in $W$
which  it represents; we rely on the context to clarify our intentions.

Given Coxeter generators $s, s^\prime \in S$, we let $m(s,s^\prime)$
denote the (possibly infinite) order of $ss^\prime$. This information
is gathered in the {\em Coxeter graph} whith vertex set $S$ and an
edge labelled $m(s,s^\prime)$ connecting $s$ and $s^\prime$ if and only if they
do not commute (i.e.\ if and only if $m(s,s^\prime)\geq 3$). As is
customary, omitted edge labels are understood to equal $3$.

The set of {\em reflections} in $W$ is $T=\{w^{-1}sw\mid s\in
S\}$. The {\em absolute length} $\ell^\prime(w)$ is the minimal $k$
such that $w=t_1\cdots t_k$ for some $t_i\in T$.

Given $w\in W$, we define its {\em (right) descent set} by 
\[
D_R(w) = \{s\in S\mid \ell(ws)<\ell(w)\}.
\]

For $J\subseteq S$, let $W_J$ be the standard parabolic subgroup
of $W$ generated by $J$. If $W_J$ is finite, it has a longest element
which is denoted by $w_J$. In $W_J$, it is characterised by the fact that
$D_R(w_J)=J$. Following tradition, we write $w_0 = w_S$.

We now define the Bruhat order. Among
several equivalent definitions, the one which follows is probably best
suited for our purposes. A disadvantage is that it is not obvious that
what it defines is indeed a partial order.

\begin{definition} \label{de:bruhat}
The {\em Bruhat order} on $W$ is the partial order defined by $u\leq
v$ if and only if for every reduced expression $s_1\cdots
s_k$ for $v$ there exist $1\leq i_1<i_2<\cdots<i_m\leq k$ such that
$u = s_{i_1}\cdots s_{i_m}$.
\end{definition}

\subsection{Twisted involutions in Coxeter groups}
By an automorphism of $(W,S)$ we mean a group automorphism of
$W$ which preserves $S$. Such an automorphism is determined by the
corresponding automorphism of the Coxeter graph. 

Let $\theta$ be an involutive automorphism of $(W,S)$. Recall from the
introduction that the set of {\em twisted involutions} in $W$ is
\[
\twist = \{w\in W\mid\theta(w)=w^{-1}\}.
\]
Observe that $\mathfrak{I}(\id)$ is the set of ordinary involutions.

The combinatorics of $\twist$ can be described in terms of ``reduced
expressions'' in a way which is remarkably similar to that of $W$
itself. We now proceed to review the parts of this theory that we will
use in the sequel. Most of this appeared in \cite{hultman2}, but
some results are new.

Define a set of symbols $\SSS=\{\s \mid s\in S\}$. There is an action
of the free monoid $\SSS^*$ on the set $W$ defined by 
\[
w\s = 
\begin{cases}
ws & \text{ if }\theta(s)ws = w,\\
\theta(s)ws & \text{ otherwise},
\end{cases}
\]
and $w\s_1\cdots\s_k = (\cdots((w\s_1)\s_2)\cdots)\s_k$. Abusing notation, we
will write $\s_1\cdots \s_k$ for $e\s_1\cdots \s_k$, where $e\in
W$ is the identity. The relevance of all this is the following:

\begin{proposition}[Proposition 3.5 in \cite{hultman2}]
The orbit of $e$ under the $\SSS^*$-action is $\twist$, i.e.\ the
twisted involutions are the elements of the form $\s_1\cdots \s_k$.
\end{proposition}

For $w\in \twist$, the {\em rank} $\rho(w)$ is the smallest $k$ such that
$w=\s_1\cdots \s_k$ for some $s_i\in S$. Then, the expression $\s_1\cdots
\s_k$ is called a {\em reduced $\SSS$-expression} for $w$.

It follows from \cite[Theorem 4.8]{hultman1} that $\Br(\twist)$, the
Bruhat order on twisted involutions, is graded with rank function
$\rho$. More precisely, it was shown that the rank function is
$(\ell+\ell^\theta)/2$, where $\ell^\theta$ is the {\em twisted
  absolute length}, but it follows from the proof that this function
coincides with $\rho$. The twisted absolute length was defined
differently in \cite{hultman1}, but here is a description which is
more convenient for our purposes:

\begin{proposition} \label{pr:elltwist}
Suppose $\s_1\cdots \s_k$ is a reduced $\SSS$-expression for $w\in
\twist$. Then, $\ell^\theta(w)$ is the number of indices $i\in [k] =
\{1, \dots, k\}$ that satisfy $\s_1\cdots \s_{i-1}\s_i$ = $\s_1\cdots
\s_{i-1}s_i$. 
\begin{proof}
Let $\lambda(w)$ be the asserted number (a priori, it depends on the
choice of reduced $\SSS$-expression). One can check that
$\ell(\theta(s)ws) = \ell(w) \Leftrightarrow \theta(s)ws = w$ for all $w\in
\twist$ and $s\in S$ (see \cite[Lemma 3.4]{hultman2}). Furthermore, by
\cite[Lemma 3.8]{hultman2},
$\rho(w\s) > \rho(w) \Leftrightarrow \ell(ws) > \ell(w)$. The
construction of the $\SSS^*$-action therefore implies    
\[
\ell(w) = 2\rho(w)-\lambda(w). 
\]
Thus, $\lambda + \ell = 2\rho = \ell^\theta + \ell$.
\end{proof}
\end{proposition}

As the terminology suggests, $\ell^\id =\ell^\prime$; see \cite{hultman1}. 

Either by a simple induction argument based on Proposition
\ref{pr:elltwist} or as an immediate consequence of \cite[Definition
  4.5]{hultman1}, we conclude that the twisted identities are the twisted
involutions of vanishing twisted absolute length:
\[
\iot = \{w\in \twist\mid \ell^\theta(w)=0\} = \{w\in \twist\mid
\ell(w) = 2\rho(w)\}.
\] 
A useful consequence is that $w\s \in \iot$ whenever $w\in \iot$ and
$s\in D_R(w)$.

Twisted involutions have reduced $\SSS$-expressions of a convenient
form, as shown by the next lemma which is due to Springer \cite{springer}. 

\begin{lemma}[Proposition 3.3(a) in \cite{springer}] \label{le:springer}
Any $u\in \twist$ with $\ell^\theta(u)=k$ has a reduced
$\SSS$-expression of the form $\s_1\cdots \s_t$ with the property that
for some $0\leq i\leq t$, $\s_1\cdots \s_i$ is a reduced
$\SSS$-expression for the longest element $w_J$ in a
$\theta$-invariant standard parabolic subgroup $W_J$, $J\subseteq S$,
with $\ell^\theta(w_J)=k$.    
\end{lemma}

The following fundamental lemma will be put to use repeatedly
throughout the paper. It is completely analogous to the corresponding
property in $W$ (due to Deodhar \cite[Theorem 1.1]{deodhar}) and
explains why $\Br(\twist)$ and $\Br(W)$ are so similarly behaved. The
first two parts are \cite[Lemma 3.9]{hultman2}. Together, they imply
the third.

\begin{lemma}[Lifting property] 
Let $s\in S$ and $u,w\in \twist$ with $u\leq w$. Suppose $s\in D_R(w)$. Then,
\begin{enumerate}
\item[(i)] $u\s \leq w$.
\item[(ii)] $s \in D_R(u) \Rightarrow u\s \leq w\s$.
\item[(iii)] $s\not \in D_R(u) \Rightarrow u \leq w\s$.
\end{enumerate}
\end{lemma}

Next, we show that $\Br(\twist)$, like $\Br(W)$, can be described in
terms of subwords. The result extends Richardson and Springer's
\cite[Corollary 8.10]{RS1}\footnote{More precisely, \cite[Corollary
    8.10]{RS1} applies to a partial order on $\twist$ which later
was shown to coincide with $\Br(\twist)$ in \cite{RS2}.} from finite to
general Coxeter groups (and adjusts it to our definition of the
$\SSS^*$-action which differs somewhat from the one given in \cite{RS1}). 

\begin{theorem}[Subword property for $\twist$]\label{th:subword}
Let $u,v\in \twist$ be given. Then, $u\leq v$ in the Bruhat order if and only if
for every reduced $\SSS$-expression $\s_1\cdots \s_k$ for $v$ there
exist $1\leq i_1 < i_2 < \cdots < i_m \leq k$ such that $u = \s_{i_1}
\cdots \s_{i_m}$.
\begin{proof}
We begin with the ``only if'' direction, so suppose $u \leq v$ and choose
a reduced $\SSS$-expression $\s_1\cdots \s_k$ for $v$. If
$s_k \not \in D_R(u)$, then the lifting property shows that $u \leq
v\s_k = \s_1\cdots \s_{k-1}$. By induction on the rank of $v$, we are
done. Otherwise, if $s_k \in D_R(u)$, we have $u\s_k \leq
v\s_k$. It follows, again by induction, that $u\s_k = \s_{i_1} \cdots
\s_{i_m}$ for some $1\leq i_1<i_2<\cdots<i_m\leq k-1$. Acting by $\s_k$
yields the desired result. 

Now consider the ``if'' part of the assertion. Assume $\s_1\cdots \s_k$
is a reduced $\SSS$-expression for $v$ and let $u = \s_{i_1} \cdots
\s_{i_m}$ for some $1\leq i_1<i_2<\cdots<i_m\leq k$.

Given $t\in [k]$, let $v_t = \s_1\cdots \s_t$ and $u_t = \s_{i_1} \cdots
\s_{i_l}$, where $i_l\leq t < i_{l+1}$ (with $i_0=1$,
$i_{m+1}=k+1$). Assuming by induction that $u_{t-1} \leq v_{t-1}$, we
either have $u_t = u_{t-1} \leq v_{t-1} < v_{t-1}\s_t=v_t$ or (by the lifting
property) $u_t = u_{t-1}s_t \leq v_{t-1}\s_t = v_t$. Thus, $u = u_k
\leq v_k = v$.   
\end{proof} 
\end{theorem}

\section{Manifestations of $\iot$}\label{se:motivation}
Henceforth, let $(W,S)$ be a finitely generated Coxeter system with
an involutive automorphism $\theta$. Recall that $\theta$ is determined by
an automorphism of the Coxeter graph. 

It should be clear by now that the goal of this paper is to study
certain properties of $\iot$. In this section, we motivate this task
by indicating some situations in which $\iot$ naturally shows up.

\begin{figure}[htb]
\epsfig{height=9cm, width=8cm, file=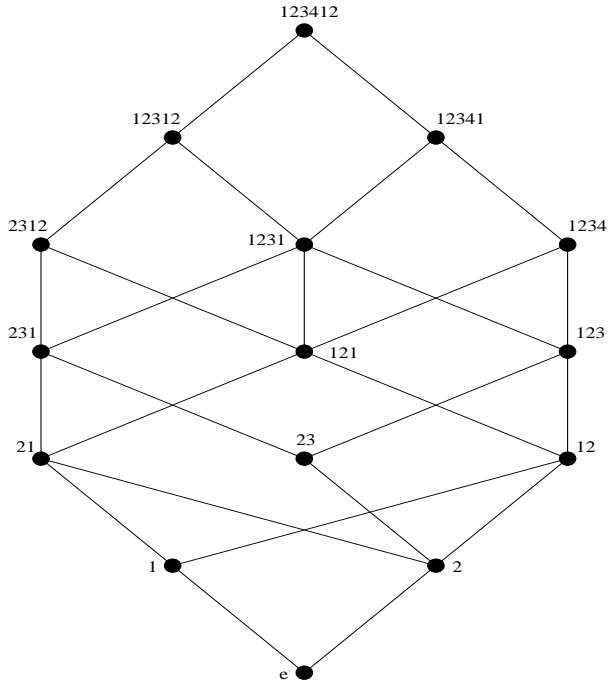}
\caption{A picture of the poset $\Br(\iot)$ when $W$ is the symmetric group
  $\mathfrak{S}_6 \cong A_5$ and $\theta$ is the unique
  non-trivial automorphism (which sends the transposition $s_i
  = (i,i+1)$ to $s_{6-i}$). The vertex labels in the figure are index
  sequences of reduced 
$\SSS$-expressions for the corresponding twisted identities. For
example, $2312$ indicates the twisted identity $\s_2\s_3\s_1\s_2$.}
\label{fi:A5} 
\end{figure}

\begin{example}
{\em Consider a connected reductive linear
algebraic group $G$ (over $\C$, say). Let $B$ be a Borel subgroup and
$T \subseteq B$ a maximal torus. Given an involutive automorphism
$\Theta:G\to G$ which preserves $T$ and $B$, let $K$ be the fixed point
subgroup. By means of left translations, $B$ acts on the symmetric
variety $G/K$. This gives rise to a finite number of orbits that may
be ordered by containment of their Zariski closures. Richardson and
Springer \cite{RS1, RS2} studied this poset $V$ by defining an order preserving
map $\varphi:V \to \Br(\twist)$. Here, the underlying Coxeter group
$W$ is the Weyl group $N/T$ (where $N$ is the normaliser of $T$ in
$G$) and $\theta$ is induced by $\Theta$.

In general, $\varphi$ is neither injective nor surjective. However,
$\iot$ is always contained in the 
image. Moreover, $\varphi$ produces an isomorphism $V \cong \Br(\iot)$
in certain interesting cases. For instance, with $G =
\mathrm{SL}_{2n}$, we may define $\Theta$ so that $K \cong
\mathrm{Sp}_{2n}$ as in \cite[Example 10.4]{RS1}. The corresponding
poset $V$, governing the orbit 
  decomposition of $\mathrm{SL}_{2n}/\mathrm{Sp}_{2n}$, is then
    isomorphic to $\Br(\iot)$, where $W$ is the symmetric group
    $\mathfrak{S}_{2n}$ and $\theta$ is given by conjugation with the longest
    element in $W$ (i.e.\ the reverse permutation $i \mapsto
    2n+1-i$). Figure \ref{fi:A5} displays this poset for $n=3$.}
\end{example}

\begin{example}[cf.\ Example 10.1 in \cite{RS1} and Example 3.2 in
    \cite{hultman2}] \label{ex:cox}
{\em Suppose $\theta$ is the involution $\theta:
W\times W \to W\times W$ given by $(v,w) \mapsto (w,v)$, where $W$ is any
Coxeter group. Observe that $\iot = \{(w,w^{-1})\mid w \in
W\}$. Hence, we have a bijection $W \leftrightarrow \iot$ in this
case. Furthermore, it is clear from Definition \ref{de:bruhat} that
this bijection gives a poset isomorphism 
$\Br(\iot) \cong \Br(W)$. Therefore, (Bruhat orders on) twisted identities
generalise (Bruhat orders on) Coxeter groups.}
\end{example}

\begin{example}\label{ex:quotient}
{\em Let $\fix$ denote the fixed point
subgroup of $\theta$. It is known (\cite{hee,muhlherr,steinberg}) that
$\fix$ itself is a Coxeter group with the following canonical set of
Coxeter generators:  
\[
\{w_J\mid J=\{s,\theta(s)\},\, m(s,\theta(s))<\infty , \, s\in S\}.
\]

Observe that $\theta(w^{-1})w = \theta(v^{-1})v \Leftrightarrow 
vw^{-1} \in \fix \Leftrightarrow v \in \fix w$. In other words, there
is a bijection between $\iot$ and the set of cosets $\fix \backslash
W$. Thus, $\iot$ can be thought of as a quotient of Coxeter groups.}
\end{example}

\begin{example}\label{ex:conjugation}
{\em If $W$ is finite, it contains a unique
longest element $w_0$. Then, $\theta(x) = w_0xw_0$ defines an
involutive automorphism of $(W,S)$. In this situation, $u \in \iot
\Leftrightarrow u = w_0w^{-1}w_0w$ for some $w\in W$. This means that
$u \in \iot$ if and only if $w_0u$ and $w_0$ are conjugate. Now, (left
or right) multiplication by $w_0$ is an antiautomorphism of
$\Br(W)$. The poset $\Br(\iot)$ is therefore isomorphic to the dual of
the subposet of $\Br(W)$ induced by the conjugacy class of the longest
element. 

In the special case when $W$ is the symmetric group $\mathfrak{S}_{2n}$, the
conjugacy class of $w_0$ consists of all fixed point free involutions
of $[2n]$. Equivalently, if we think of a permutation in terms of its
disjoint cycle decomposition, the conjugacy class of $w_0$ corresponds
to the set of complete matchings on $2n$ elements. Thus, $\Br(\iot)$ in
this case gives (the dual of) a Bruhat order on matchings. 

We mention in this context that the literature already contains a
``Bruhat order'' on matchings defined by Deodhar and Srinivasan
\cite{DS}. Their poset is, however, strictly weaker than the dual of
$\Br(\iot)$, although both posets are graded with the same rank
function; see Remark \ref{rem:poin1}. }
\end{example}

\section{The Bruhat order on twisted identities}\label{se:bruhat}

We now turn to the core of the paper, namely the study of
$\Br(\iot)$. To begin with, we observe that it sometimes coincides with
the more familiar $\Br(\twist)$.

\begin{proposition}\label{pr:coincide}
If there exists no $s\in S$ such that $\theta(s)$ and $s$ are
conjugate, then $\iot = \twist$. 
\begin{proof}
Notice that $\theta(s)ws = w$ if and only if $w^{-1}\theta(s)w =
s$. Thus, if $s$ and $\theta(s)$ are not conjugate, we have $w\s =
\theta(s)ws$ for all $w$. Now, the assertion is a consequence of Proposition
\ref{pr:elltwist} and the remarks that follow it. 
\end{proof}
\end{proposition}

The hypothesis of Proposition \ref{pr:coincide} means that for
all $s\in S$, every path from $s$ to $\theta(s)$ in the Coxeter graph of $W$
contains an edge with an even label; see \cite[Exercise 1.16]{BB}. For
example, Proposition 
\ref{pr:coincide} applies when $W=F_4$ with
$\theta$ being the unique non-trivial automorphism.

We will find that $\iot$ is particularly well-behaved if $\theta$
never flips an edge with an odd label in the Coxeter graph of $W$. 

\begin{definition}
Say that $\theta$ has the {\em no odd flip (NOF) property} if
$m(s,\theta(s))$ is even or infinite for all $s \neq \theta(s)\in S$.
\end{definition} 

For example, if $W$ is finite and irreducible, $\theta$ has the NOF
property unless $W$ is of type $A_{2n}$ or $I_2(2n+1)$, $n\in
\N$, and $\theta$ is the unique non-trivial automorphism.

We now proceed to prove a series of lemmata that put restrictions on the
behaviour of elements of twisted absolute length $1$ when $\theta$ has
the NOF property. Informally put, these restrictions
ensure that intervals in $\Br(\iot)$ never can be ``sparse enough''
not to inherit the gradedness from $\Br(\twist)$.

\begin{lemma} \label{le:length1}
Suppose $W$ is irreducible and finite. If $\ell^\theta(w_0)=1$, then
$W$ is either the rank one group $A_1$ or the dihedral group $I_2(m)$
for some odd integer $m$. 
\begin{proof}
If $\theta=\id$, the twisted absolute length coincides with the absolute
length. Hence, $w_0$ is a reflection. It follows, for example by the
classification of finite Coxeter groups, that $W=A_1$ or
$W=I_2(m)$, $m$ odd.

Now, assume $\theta\neq \id$. Obviously, $\ell^\theta(w)$ and
$\ell(w)$ always have the same parity. In particular, $\ell(w_0)$ is
odd. By inspection of Coxeter graphs, the only groups with $\ell(w_0)$
odd that admit a nontrivial
$\theta$ are $I_2(m)$, $m$ odd, and $A_m$, $m \equiv 1,2
\pmod{4}$. These groups have $\theta(w) = w_0ww_0$ as unique nontrivial
automorphism; see \cite[Exercise 4.10]{BB}. Therefore, $w \mapsto
w_0w$ is a bijection which sends $\Br(\twist)$ to the dual of
$\Br(\twistid)$. In particular, the top
element $w_0$ has
the same rank in $\Br(\twist)$ as in $\Br(\twistid)$. This implies
$\ell^\prime(w_0) = \ell^\theta(w_0) = 1$, i.e.\ that $w_0$ is again a
reflection. Hence, $W$ is the dihedral group $I_2(m)$ for some odd $m$.
\end{proof}
\end{lemma}

\begin{lemma} \label{le:length2}
Let $W$ be finite and assume $\theta$ has the NOF property. If
$\ell^\theta(w_0)=1$, then $w_0$ has a reduced 
$\SSS$-expression beginning with $\s$ for some $\theta$-fixed $s\in S$.
\begin{proof}
The twisted absolute length clearly is additive over direct products
of $\theta$-invariant Coxeter systems. Thus, $W$ decomposes as $W = W_J \times
W_{S\setminus J}$ where $\ell^\theta(w_J)=1$,
$\ell^\theta(w_{S\setminus J})=0$, $W_J$ is irreducible and
$J\subseteq S$ is invariant under $\theta$. Choose a reduced
$\SSS$-expression $\s_1\cdots \s_k$ for $w_{S\setminus J}$.

Lemma \ref{le:length1} shows that either
$W_J \cong A_1$ or $W_J\cong I_2(m)$, $m$ odd. In the former case, $w_J =
\s$ for some $\theta$-fixed $s\in S$. This implies that $w_0$ has the
reduced $\SSS$-expression $\s\s_1\cdots\s_k$. In the latter
situation, $w_J$ has the reduced $\SSS$-expression
$\s\s^\prime \s\s^\prime \cdots$ ($(m+1)/2$ letters), where
$J=\{s,s^\prime\}$. The hypothesis on 
$m(s,\theta(s))$ means 
that $s$ and $s^\prime$ are $\theta$-fixed. Hence,
$\s\s^\prime \s\s^\prime \cdots \s_1\cdots\s_k$ is a reduced expression for
$w_0$ of the desired form.
\end{proof}
\end{lemma}

\begin{lemma} \label{le:cover}
Suppose $\theta$ has the NOF property. Then, given $v\in \twist$ with
  $\ell^\theta(v)=1$, $v$ covers at most one twisted identity in $\Br(\twist)$.
\begin{proof}
By Lemma \ref{le:springer}, $\ell^\theta(v)=1$ implies $v = w_J
\s_1\cdots\s_k$ with $\ell(v)=\ell(w_J)+2k$ for some $\theta$-invariant
$J\subseteq S$ satisfying $\ell^\theta(w_J)=1$. Invoking Lemma
\ref{le:length2}, we may assume $J=\{s\}$. In other words, $v$ has a reduced
$\SSS$-expression of the form $\s \s_1\cdots \s_k$, where $s$ is fixed
by $\theta$. 

Each reduced $k$-letter subword of $\s\s_1\cdots \s_k$ except at most
one begins with $\s$, hence represents an element of non-vanishing twisted
absolute length. The exception is $\s_1\cdots \s_k$ (if this expression
is reduced) which, by Theorem \ref{th:subword}, is the only candidate for
a twisted identity covered by $v$. 
\end{proof}
\end{lemma}

Next, we state our main result on gradedness. As a notable special
case, $\Br(\iot)$ is graded whenever the Coxeter graph of $W$ is a
tree containing a $\theta$-fixed node. The nature of the rank function
should not be surprising. Indeed, if $v=\s_1\cdots \s_k$ with
$\rho(v)=k$, then $e<\s_1<\s_1\s_2 < \cdots < \s_1\cdots \s_k$ is an
unrefinable chain in $\Br(\iot)$. Hence, the rank function necessarily
is $\rho$ whenever $\Br(\iot)$ is graded.

\begin{theorem} \label{th:graded} If $\theta$ has the NOF property,
  then $\Br(\iot)$ is graded with rank function $\rho$.
\begin{proof}
Since we know that $\Br(\twist)$ is graded with rank function $\rho$,
it is enough to show that every interval $[u,v]\subseteq \Br(\iot)$
contains a chain of length $\rho(v)-\rho(u)$.

In order to get a contradiction, assume the theorem is false. Choose
a minimal interval $[u,v] \subseteq \Br(\iot)$ which does not contain
a chain of length $\rho(v)-\rho(u)$ and such that $\rho(v)$ is minimal
among all such intervals. 

Pick $s\in D_R(v)$. We must have $v\s \not \geq u$; otherwise $[u,v]$
would not be minimal. By the lifting property, $u\s < u,
v\s$. Minimality of $\rho(v)$ shows that $[u\s,v\s]$ contains some $a\in
\iot$ which covers $u\s$ in $\Br(\twist)$. Using the lifting property
once again, we conclude that $a\s$ 
covers $u$ in $\Br(\twist)$. Minimality of $[u,v]$ implies $a\s \not
\in \iot$. Thus, $a\s = as$ and $\ell^\theta(a\s) = \ell^\theta(a)+1
= 1$. Now, $a\s$ covers two distinct twisted identities, namely $a$
and $u$. This contradicts Lemma \ref{le:cover}, and the proof is complete.   
\end{proof}
\end{theorem}

Unfortunately, Theorem \ref{th:graded} does not give the full
picture. For example, one readily checks that with $W$ being any
dihedral group and $\theta$ the non-trivial automorphism, $\Br(\iot)$
is graded, although
this is not predicted by Theorem \ref{th:graded}. At this
point one may very well suspect that $\Br(\iot)$ is always graded. This is
not true, as shown by the next example.

\begin{example}\label{ex:ungraded} {\em 
Suppose $W\cong \widetilde{A_2}$, the affine Weyl group corresponding
to $A_2$. Its Coxeter generating set is $S=\{s_1,s_2,s_3\}$, and
$m(s_1,s_2)=m(s_1,s_3)=m(s_2,s_3)=3$. Define the involution $\theta$
by $\theta(s_1)=s_1$, $\theta(s_2)=s_3$ and $\theta(s_3)=s_2$.

Consider the twisted identities $u=\s_3=s_2s_3$ and
$v=\s_2\s_1\s_3=s_2s_1s_3s_2s_1s_3$. As an interval in $\Br(\twist)$,
the open interval $(u,v)$ contains two elements, namely $\s_2\s_3 =
s_3s_2s_3$ and $\s_1\s_3=s_2s_1s_3$, neither being a twisted
identity. Thus, $e<u<v$ and $e<\s_2<\s_2\s_1<v$ are
unrefinable chains in $\Br(\iot)$. Hence, it is not graded.}
\end{example}

Our attention is now turned to the topology of intervals in
$\Br(\iot)$. The first results are valid without any restriction on $W$ or
$\theta$.

\begin{theorem}\label{th:collapsible}
Consider an interval $[u,v] \subseteq \Br(\iot)$. If $u\s = us$
for some $s\in D_R(v)$, then $(u,v)$ is collapsible.
\begin{proof}
Let $I = (u,v)$. Our strategy is to invoke Theorem \ref{th:forman} after
defining a complete acyclic matching $M$ on the face poset $P=P(\Delta(I))$. 

Suppose $c$ is a chain in $I$, and define $x_c = \min \{x\in
c\cup\{v\}\mid s\in D_R(x)\}$. Now, let $M:P\to P$ be given by
\[
M(c) = 
\begin{cases}
c \cup \{x_c\s\} & \text{ if }x_c\s \not \in c,\\
c\setminus \{x_c\s\} & \text{ otherwise.}
\end{cases}
\]
Observe that $x_c\s \neq u$ because $u\s = us \not \in \iot$. The
lifting property ensures that $x_c\s \cup c$ is a chain in $I$. Since
$s\not \in D_R(x_c\s)$, we have $x_c = x_{M(c)}$ and therefore
$M(M(c))=c$. Hence, $M$ is a complete matching.

To prove acyclicity of $M$, suppose we have four chains $c \prec
M(c) \succ c^\prime \prec M(c^\prime)$, where $\prec$ is the covering
relation in 
$P$ and $c\neq c^\prime$. By construction, $M(c) = c\cup \{x_c\s \}$
and $c^\prime =M(c)\setminus \{x\}$ for some $x\in M(c)$. In fact, $x
= x_c = x_{M(c)}$; otherwise we would have $x_{c^\prime} = x_c$,
implying $M(c^\prime) = c^\prime \setminus \{x_c\s\} \prec
c^\prime$. We therefore conclude that $c^\prime = \{x_c\s\}\cup c 
  \setminus \{x_c\}$. In particular, $c^\prime$ has strictly
  fewer elements with $s$ as a descent than $c$ does. This shows the
  impossibility of a sequence
\[
c_1 \prec M(c_1)\succ c_2\prec M(c_2)\succ \cdots \prec M(c_{t-1})\succ
c_t = c_1
\]
with $c_1\neq c_2$. In other words, $M$ is acyclic.
\end{proof}
\end{theorem}

The strategy employed in the preceding proof is not quite applicable if $u\s
= \theta(s)us$, $s\not \in D_R(u)$, because chains that contain $u\s$
would not have a well-defined matching partner. Leaving 
these chains unmatched, however, allows us to conclude a result which
is useful in recursive arguments.

\begin{theorem}\label{th:recursive}
Consider an interval $[u,v] \subseteq \Br(\iot)$. If $u\s =
\theta(s)us$ for some $s \in D_R(v)\setminus D_R(u)$, then $(u,v)$
is homotopy equivalent to the suspension of $(u\s,v)$.
\begin{proof}
Just as in the proof of Theorem \ref{th:collapsible}, we define a
matching on $P=P(\Delta(I))$, where $I=(u,v)$. This time, there will
be critical cells. Define 
\[
Z = \{c\in \Delta(I)\mid c \ni u\s\}.
\]
Now, arguing exactly as in the proof
of Theorem \ref{th:collapsible} we have a complete acyclic matching
$M:P\setminus Z \to P \setminus Z$.

The set of critical cells is the order filter $Z$ which consists of
the cells in the open star of $u\s$. Corollary \ref{co:filter}
therefore implies that $\Delta(I)$ is homotopy equivalent to the star
of $u\s$ with the entire link of $u\s$ collapsed to a single point. The complex
obtained in this way is homotopic to the suspension of the link of
$u\s$, and this link is $\Delta((u\s,v))$.  
\end{proof} 
\end{theorem}

Using Theorem \ref{th:recursive}, we can always compute the homotopy
type of an interval in terms of that of another interval. The problem
is that if $D_R(v)\subseteq D_R(u)$, we are
only able to express $(u,v)$ in terms of larger intervals. This
makes it difficult to set up an 
inductive argument since we may not reach a reasonable base case. To
rectify the situation, we pay the price of restricting $W$ and $\theta$.

Let us call an interval $[u,v]\subseteq \Br(\iot)$ {\em full} if
$\{w\in \twist\setminus \iot \mid u\leq w \leq v\} = \emptyset$. 

\begin{lemma} \label{le:full}
Suppose $\theta$ has the NOF property. If an interval $[u,v]\subseteq
\Br(\iot)$ is full and $s\not \in 
D_R(v)\cup D_R(u)$, then either $u\s = us$ or $[u,v\s]$ is full
(implying, in particular, $v\s \in \iot$).
\begin{proof}
Assume, in order to deduce a contradiction, that $u\s = \theta(s)us$ and
$[u,v\s]$ is not full. Choose $w\in [u,v]$ minimal with the property
that $w\s = ws$; it exists since every element in $[u,v\s]\setminus
[u,v]$ is of the form $w\s$, $w\in [u,v]$. We have $\ell^\theta(ws)=1$
and $u< w < v\s$. 

Since $u< w< ws$ and intervals in $\Br(\twist)$ are PL spheres, $ws$
must cover some element $w^\prime\in \twist$ with $u<w^\prime\neq
w$. By the lifting property, $u\leq w^\prime\s < w$. Since $[u,v]$ is
full, $w^\prime\s \in \iot$. Minimality of $w$ now implies
$w^\prime \in \iot$. Hence, $ws$ covers
more than one twisted identity, contradicting Lemma \ref{le:cover}.  
\end{proof}
\end{lemma}

\begin{lemma}\label{le:matching}
Suppose $\theta$ satisfies the NOF property. Let $u<v$ for $u,v\in
\iot$ and suppose $v\s = vs$ for some $s\in 
S$. If $u\s = \theta(s)us > u$, then $u\s < v$.
\begin{proof}
We employ induction on $r = \rho(v)-\rho(u)$. The assertion is vacuously
true if $r = 1$, because in this case the hypotheses
imply that $vs$ covers both $u\s$ and $v$, contradicting Lemma \ref{le:cover}.

Now suppose $r>1$. Since $\Br(\iot)$ is graded, we can choose $w\in
\iot$ such that $w$ covers $u$ and $w<v$. If $u\s = w$, we are
done. Otherwise, the lifting property implies that $w\s$ covers
$u\s$. As in the $r=1$ case above, $w\s = ws$ is impossible, so $w\s =
\theta(s)ws$. By the induction assumption, $w\s < v$. Hence, $u\s < v$.
\end{proof}
\end{lemma}

\begin{theorem}\label{th:topology}
Assume $\theta$ has the NOF property. Consider an interval
$[u,v]\subseteq \Br(\iot)$. If it is full, 
then $(u,v)$ is a PL sphere of dimension
$\rho(v)-\rho(u)-2$. Otherwise, $(u,v)$ is $\Z$-acyclic.
\begin{proof}
If $[u,v]$ is full, the assertion follows immediately from the
corresponding result on $\Br(\twist)$, namely \cite[Corollary
  4.6]{hultman2}.

Suppose $[u,v]$ is not full, and let $I=(u,v)$. Proceeding by induction
on $\rho(v)$, we choose $s\in D_R(v)$. There are two cases depending
on whether or not $s$ also is a descent of $u$.

{\bf Case I, $s\not \in D_R(u)$.}

If $u\s = us$, Theorem \ref{th:collapsible} shows that $I$ is
 collapsible, in particular contractible and $\Z$-acyclic. By Lemma
 \ref{le:full}, 
 this is always the case if $[u,v\s]$ is full. Therefore, we may
 assume that $[u,v\s]$ is not full and $u\s = \theta(s)us$. By the
 induction assumption, $(u,v\s)$ is $\Z$-acyclic.

Let $P=P(\Delta(I))$. Given a chain $c \in P$, define $x_c =
\max\{x\in c\cup\{u\}\mid x\s = 
\theta(s)xs > x\}$; this set is nonempty since it contains $u$. Let
$Z = \{c\in P\mid v\s \in c\} = \{c\in P\mid x_c = v\s\}$ and
define $M:P\setminus Z\to P\setminus Z$ by
\[
M(c) = 
\begin{cases}
c \cup\{x_c\s\} & \text{if } x_c\s \not \in c,\\
c \setminus \{x_c\s\} & \text{if } x_c\s \in c.
\end{cases}
\]
Lemma \ref{le:matching} together with the lifting property
proves that $c\cup \{x_c\s\} \in P$. Observing that $x_c = x_{M(c)}$,
we conclude that $M$ is a matching on $P$ with set of critical cells
$Z$. An argument completely analogous to the acyclicity part of the
proof of Theorem \ref{th:collapsible} shows that $M$ is
an acyclic matching. Exactly as in the proof of Theorem
\ref{th:recursive}, this implies that $I$ is homotopy equivalent to
the suspension of the $\Z$-acyclic complex $(u,v\s)$. Since a complex
and its suspension have isomorphic reduced homology groups up to an
index shift, $I$ is $\Z$-acyclic.

{\bf Case II, $s\in D_R(u)$.}

By Lemma \ref{le:full}, $[u\s,v\s]$ is not full, and therefore
$(u\s,v\s)$ is $\Z$-acyclic. Combining Case I with Theorem
\ref{th:recursive}, we have the following homotopy equivalences:
\[
\Sigma I \simeq (u\s,v) \simeq \Sigma (u\s,v\s),
\]
where $\Sigma$ denotes suspension. Hence, the suspension of $I$ is
$\Z$-acyclic. Therefore, $I$ is $\Z$-acyclic, too.
\end{proof}
\end{theorem}

\section{The Poincar\'e series of $\iot$} \label{se:poincare}

From a combinatorial point of view, the {\em Poincar\'e series} of $W$
is simply the length generating function:
\[
\poi(W;t) = \sum_{w\in W}t^{\ell(w)}.
\]
Analogously, we may define the Poincar\'e series of $\iot$ to be
\[
\poi(\iot;t) = \sum_{w\in \iot}t^{\rho(w)}.
\]
In particular, $\poi(\iot;t)$ is the rank generating function for
$\Br(\iot)$ whenever it is graded.

Since $\iot$ is, in some sense, a quotient of Coxeter groups (cf.\
Example \ref{ex:quotient}), one may naively hope that
its Poincar\'e series is a quotient of the series of the
corresponding groups; let us say that $\poi(\iot;t)$ {\em
  factors through $\theta$} if $\poi(W;t)=\poi(\iot;t)\poi(\fix;t)$. 

With $W=A_2$ and $\theta$ given by conjugation with the
longest element, for example, one obtains $\poi(\iot;t)=1+2t$, whereas
$\poi(W;t)/\poi(\fix;t)=1+t+t^2$. Thus, in this case, $\poi(\iot;t)$
does not factor through $\theta$. Intriguingly, though, factorisation
does occur in several situations, as we shall see below. 

First, we give a necessary condition for the Poincar\'e series of
$\iot$ to factor through $\theta$. The condition is that
$\theta$ (seen as a Coxeter graph automorphism) is not allowed to flip
edges unless they are labelled $\infty$. 

\begin{proposition}\label{pr:factor}
If $\poi(W;t)=\poi(\iot;t)\poi(\fix;t)$, then $m(s,\theta(s))\in
\{1,2,\infty\}$ for all $s\in S$.
\begin{proof}
Suppose equality holds, and consider the coefficient of the linear
term on both sides. On the left hand side, this coefficient is
$|S|$. The elements of length $1$ in $\fix$ (i.e.\ the Coxeter
generators) are the longest elements in the finite parabolic subgroups
of the form $\langle s,\theta(s)\rangle$. In $\iot$, the elements of
rank $1$ are the twisted identities of the form $\s$, $s\in S$. Now,
$\s$ is a twisted identity if and only if $\theta(s)\neq s$. Moreover, for $s\neq
s^\prime$, $\s =
\s^\prime$ if and only if $\theta(s) = s^\prime$ and $m(s,\theta(s))=2$. As a
consequence, the linear coefficient of the right hand side is 
\[
|\{\{s,\theta(s)\}\mid m(s,\theta(s))<\infty\}| + |S| - |\{s \mid
 s=\theta(s)\}| - |\{\{s,\theta(s)\}\mid m(s,\theta(s))=2\}|.
\]
This is equal to $|S|$ if and only if $m(s,\theta(s))\in
\{1,2,\infty\}$ for all $s\in S$.
\end{proof}
\end{proposition}

It turns out that in the finite setting, Proposition \ref{pr:factor}
actually gives a characterisation of the $(W,\theta)$ for which the
factorisation phenomenon occurs. (The label $\infty$, of course,
cannot occur in this case).

\begin{figure}[htb]
\epsfig{height=9cm, width=10cm, file=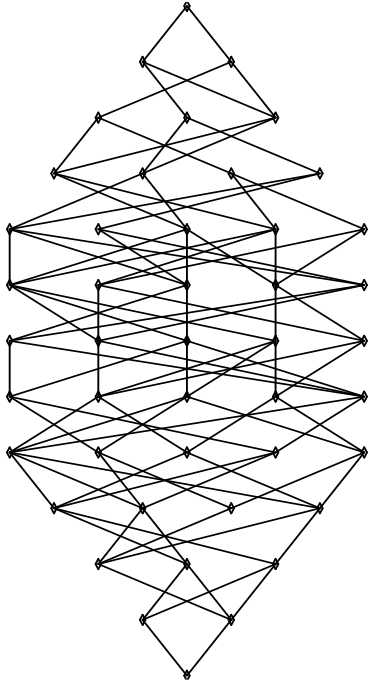}
\caption{The poset $\Br(\iot)$ when $W=E_6$ and $\theta$ is the
  non-trivial automorphism.} \label{fi:E6}
\end{figure}

\begin{theorem} \label{th:poincare}
Suppose $W$ is finite. Then, $\poi(\iot;t)$ factors through $\theta$ if and
only if $m(s,\theta(s)) \in \{1,2\}$ for all $s\in S$, i.e.\ if and
only if $s$ commutes with $\theta(s)$ for all $s \in S$.
\begin{proof}
Suppose $J\subseteq S$ is minimal such that $\theta(J)=J$ and
$m(s,s^\prime)=2$ whenever $s\in J$, $s^\prime \not \in J$. Then,
either $W_J$ is irreducible or $W_J$ is isomorphic to the direct
product of two isomorphic irreducible Coxeter groups that are
interchanged by $\theta$ as in Example \ref{ex:cox}. To prove the
theorem it therefore suffices to check that it holds in the setting of Example
\ref{ex:cox} and for the finite irreducible groups that satisfy the
hypotheses. 

To begin with, return to the setup in Example \ref{ex:cox}. Observe that $\fix
\cong W$. Furthermore, the rank of $(w,w^{-1})\in \iot$ is $\ell(w)$,
so that $\poi(\iot;t) = \poi(W;t)$. Since
$\poi(W\times W;t)=\poi(W;t)^2$, we have
$\poi(W\times W;t)=\poi(\iot;t)\poi(\fix;t)$.

It remains to consider the finite irreducible groups. If $\theta =
\id$, the assertion is trivial. If not, the groups that satisfy the
criteria are $A_{2n-1}\cong \mathfrak{S}_{2n}$, $D_n$ and $E_6$.

{\bf Type A.} Consider the symmetric group case
  $W=\mathfrak{S}_{2n}$ with 
  $\theta$ given by conjugation with the 
  longest element $w_0$ (the reverse permutation). Let us verify that
  $\poi(\iot;t)$ factors through $\theta$. 

It is known \cite[Supplementary problem 1.24]{stanley} that the 
  fixed point free involutions $F(2n)$ in $\mathfrak{S}_{2n}$ have the following
  generating function: 
\[
I(n;t) = \sum_{\pi\in F(2n)}t^{\inv(\pi)} =
t^n\prod_{i=0}^{n-1}(1+t^2+t^4+\cdots+t^{4i}), 
\]
where $\inv(\pi)$ denotes the number of inversions in $\pi$, which is
the same as the Coxeter length of $\pi$. Since $\iot = w_0F(2n)$
(Example \ref{ex:conjugation}), $\ell(w_0w) = 2n^2-n-\ell(w)$ for all
$w\in \mathfrak{S}_{2n}$ and $\rho(w) = \ell(w)/2$ for all $w\in \iot$, we
obtain
\[
\begin{split}
\poi(\iot;t) & = (t^{1/2})^{2n^2-n}I(n;t^{-1/2}) =
t^{n(n-1)}\prod_{i=0}^{n-1}(1+t^{-1}+t^{-2}+\cdots+t^{-2i})  \\
&= \prod_{i=0}^{n-1}(1+t+t^2+\cdots+t^{2i}).
\end{split}
\]
It is well-known, and straightforward to check, that $\fix \cong B_n$,
the hyperoctahedral group of rank $n$. Now,
\[
\poi(\mathfrak{S}_{2n};t) = \prod_{i=0}^{2n-1}(1+t+t^2+\cdots+t^i)
\]
and
\[
\poi(B_n;t) = \prod_{i=1}^n(1+t+t^2+\cdots+t^{2i-1}).
\]
Thus, $\poi(W;t) = \poi(\iot;t)\poi(\fix;t)$ in this situation.

{\bf Type D.} Let $W=D_n$ with the Coxeter generators being labelled
as described in the Coxeter graph below.
\begin{center}
\begin{picture}(150,60)(0,0)
\put(10,50){\circle*{3}}
\put(-1,50){$s_1$}
\put(10,10){\circle*{3}}
\put(0,4){$s_2$}
\put(30,30){\circle*{3}}
\put(30,34){$s_3$}
\put(60,30){\circle*{3}}
\put(56,34){$s_4$}
\put(130,30){\circle*{3}}
\put(126,34){$s_n$}
\put(10,50){\line(1,-1){20}}
\put(10,10){\line(1,1){20}}
\put(30,30){\line(1,0){30}}
\put(60,30){\line(1,0){10}}
\put(130,30){\line(-1,0){10}}
\put(90,29){$\dots$}
\end{picture}
\end{center}
Define $\theta$ by $\theta(s_1)=s_2$, $\theta(s_2) = s_1$ and
$\theta(s_i)=s_i$ for all $i>2$.

It is not hard to realise that $\iot$ consists precisely of the elements
of the form $\s_2\s_3\cdots \s_k$ for $k\in [n]$ (if $k=1$, we
interpret this as the identity element $e$). In particular,
$\poi(\iot;t)=1+t+t^2+\cdots + t^{n-1}$. Noting that $\fix \cong
B_{n-1}$ (whose Poincar\'e series was presented in the type $A$ case above) and
\[
\poi(D_n;t) = (1+t+t^2+\cdots+t^{n-1})\poi(B_{n-1};t),
\]
we again conclude that $\poi(\iot;t)$ factors through $\theta$.

{\bf Type E.} When $\theta$ is the unique non-trivial
automorphism of the Coxeter graph of $E_6$, the Hasse diagram of
$\Br(\iot)$ is displayed in Figure \ref{fi:E6}. Inspecting it, one
obtains the rank generating function and verifies that $\poi(E_6;t)$
factors through $\theta$. Here, $\fix \cong F_4$.
\end{proof}
\end{theorem}

\begin{remark}\label{rem:poin1}
{\em Return to the type A case $W=\mathfrak{S}_{2n}$ with 
  $\theta$ given by conjugation with $w_0$. Then, $\poi(\iot;t)$
  coincides with the rank generating function for
Deodhar and Srinivasan's ``Bruhat order'' on $F(2n)$ \cite[Theorem 1.3(i)]{DS}
which was mentioned in Example \ref{ex:conjugation}. In fact, although
formulated differently in \cite{DS}, the definition of their rank function
$\mathrm{wt}:F(2n)\to \N$ is readily seen to yield
\[
\mathrm{wt}(\sigma)= \frac{\inv(\sigma)-n}{2}=n(n-1)-\inv(w_0 \sigma)/2
= n(n-1) -\rho(w_0 \sigma),
\]
so that
\[
\sum_{\sigma\in F(2n)}t^{\mathrm{wt}(\sigma)} = t^{n(n-1)}\poi(\iot;1/t)
= \poi(\iot;t),
\]
where the last equality is valid since $\poi(\iot;t)$ is a symmetric
polynomial in this case.}
\end{remark}

\section{Further comments and open questions}\label{se:problems}

\subsection{Gradedness}
We have already commented on the fact that, although Theorem
\ref{th:graded} does not predict it in all cases, $\Br(\iot)$ is graded
whenever $W$ is dihedral. In fact, the only instances of non-graded
$\Br(\iot)$ that we know of follow from Example \ref{ex:ungraded} and simple
variations thereof\footnote{For example, the argument is not affected if
$m(s_1,s_2) = m(s_1,s_3)$ is increased. Similarly, $m(s_2,s_3)$
  can be increased to any odd number $\geq 3$ with simple
  adjustments.}. Ultimately, of course, we would like a complete 
characterisation of the $(W,S)$ and $\theta$ for which $\Br(\iot)$ is
graded. A step in this direction would be to resolve the following two
problems.

\begin{conjecture} \label{co:symmetric}
When $W=A_{2m}\cong \mathfrak{S}_{2m+1}$, $\Br(\iot)$ is
always graded. 
\end{conjecture}
Aided by Stembridge's Maple packages {\tt coxeter}
and {\tt posets} \cite{coxeter,posets}, we have verified Conjecture
\ref{co:symmetric} for $m\leq 4$. The conjecture would follow from
an affirmative answer to

\begin{question}
Is $\Br(\iot)$ graded whenever the Coxeter graph of $W$ is a tree?
\end{question}

\subsection{Interval topology}

Consider the setup of Theorem \ref{th:topology}. By a result of
Whitehead \cite{whitehead}, a collapsible PL manifold is a PL
ball (see e.g.\ \cite[Corollary 3.28]{RoSa}). By induction on the
rank, we may assume that the proper subintervals of $(u,v)$ are PL
spheres or balls, hence that $(u,v)$ is a
PL manifold. Thus, if we were somehow able to replace ``$\Z$-acyclic'' by
``collapsible'' in the conclusion of Theorem \ref{th:topology}, we would have a
proof of the next conjecture.

\begin{conjecture}\label{co:PL}
Suppose $\theta$ has the NOF property. Consider an interval
$[u,v] \subseteq Br(\iot)$. If it is full, 
then $(u,v)$ is a PL sphere of dimension
$\rho(v)-\rho(u)-2$. Otherwise, $(u,v)$ is a PL ball of the same dimension.
\end{conjecture}

The conditions on $(W,S)$ and $\theta$ which guarantee that
$\Br(\iot)$ is graded in Theorem \ref{th:graded} are exactly those for
which Theorem \ref{th:topology} asserts that the intervals are PL spheres or
$\Z$-acyclic (namely, $\theta$ should satisfy the NOF property). This
immediately leads to the next question. 

\begin{question} \label{qu:gradtop}
Does there exist a graded interval in some $\Br(\iot)$ which is
neither a PL sphere nor $\Z$-acyclic? 
\end{question}

If Conjecture \ref{co:PL} is valid, one would of course like to replace
``$\Z$-acyclic'' by ``ball'' in Question \ref{qu:gradtop}.

Like any non-graded interval, the one in Example \ref{ex:ungraded} is
neither a PL ball nor sphere. Note, however, that it is
homotopy equivalent to a $0$-dimensional sphere.

\begin{question} \label{qu:homotopy}
Are there intervals in some $\Br(\iot)$ that are neither homotopy
equivalent to spheres nor contractible?
\end{question}

The M\"obius function of an interval coincides with the reduced Euler
characteristic of its order complex. Thus, a negative answer to Question
\ref{qu:homotopy} would imply a negative answer to the following:  

\begin{question}\label{qu:Mobius}
Is there some $\Br(\iot)$ whose M\"obius function takes values outside
$\{-1,0,1\}$? 
\end{question}

When $W=A_{2m}$, $m\leq 4$, we have used \cite{coxeter,posets}
to verify that the range of the M\"obius
function of $\Br(\iot)$ is $\{-1,0,1\}$.

\subsection{Directedness}

Whenever $W$ is finite, $\Br(\twist)$ and $\Br(W)$ contain a unique
top element $w_0$. This element is, however, usually not a twisted
identity. Indeed $\Br(\iot)$ does not always possess a unique maximal
element.

\begin{proposition}\label{pr:maximal}
If $W$ is irreducible and finite, then $\Br(\iot)$ has more than one
maximal element if and only if $W$ is of type $A_{2m}$ or $I_2(2m+1)$, for
some positive integer $m$, and $\theta$ is the unique non-trivial
automorphism.
\begin{proof}
If $\theta=\id$, $e$ is the only twisted identity. Consulting the
classification of finite irreducible Coxeter groups, we find that the
groups that admit a nontrivial $\theta$ are $A_n$, $D_n$, $E_6$, $F_4$
and $I_2(n)$. Furthermore, this $\theta$ is unique in all cases
(except $D_4$, but the various choices are then equivalent). The
dihedral groups $I_2(n)$ are easy to handle. Types $D$ and $E$ are
covered by Theorem \ref{th:poincare} whereas the type $F$ assertion
follows from Proposition \ref{pr:coincide}. It remains to understand
type $A_n$. Here, one may study
the Bruhat order on the conjugacy class of $w_0$ as described in
Example \ref{ex:conjugation}. It is straightforward to check that the
minimal elements in this class are the possible products of $\lceil n/2
\rceil$ mutually commuting Coxeter generators. If $n$ is odd, there is only
one such product; otherwise there are $n/2+1$ of them.
\end{proof}
\end{proposition}

The property of having a unique maximal element has a counterpart in
infinite groups. A poset is {\em directed} if every pair
of elements has a common upper bound. 

\begin{question}\label{qu:directed}
For which $W$, $\theta$ is $\Br(\iot)$ a directed poset?
\end{question}  

We proceed to mention some reasonably straightforward partial answers
to Question \ref{qu:directed}.  

It is known that
$\Br(W)$ always is directed. A proof using the lifting property
of $\Br(W)$ is given in \cite{BB}. Employing instead the lifting property of
$\Br(\twist)$, it follows in exactly the same way that $\Br(\twist)$
is directed, too. Thus, in addition to the cases provided by Proposition
\ref{pr:maximal}, we may immediately conclude that $\Br(\iot)$ is directed
whenever the hypothesis of Proposition \ref{pr:coincide} is
satisfied. 

Here is another situation in which directedness is relatively
effortless to establish: 

\begin{proposition}\label{pr:directed}
Assume $W$ is infinite. Suppose there is a partition $S=S_1 \uplus S_2
\uplus S_3$ with $\theta(S_1) = S_2$ and 
$S_3 = \{s\in S\mid \theta(s)=s\}$ such that the elements of $S_i$
commute pairwise for $i=1,2,3$. Then, $\Br(\iot)$ is directed. 
\begin{proof}
Recall that a {\em Coxeter element} is a product of the Coxeter
generators (in any order). Speyer \cite{speyer} has shown that
$\ell(c^k)=|S|k$ for all $k\in \N$ and an arbitrary Coxeter element
$c$ whenever $W$ is infinite. Therefore, given $v,w\in
W$, a large enough power of any 
Coxeter element is an upper bound for $v$ and $w$ in
$\Br(W)$. Observe that $c=w_{S_1}w_{S_3}w_{S_2}$ is a Coxeter element
in $W$. Moreover, $c^{2k}= \theta(c^{-k})c^k \in \iot$. Thus,
$\Br(\iot)$ is directed.
\end{proof}
\end{proposition}

\subsection{The Poincar\'e series}

Proposition \ref{pr:factor} and Theorem \ref{th:poincare} immediately
lead to the next question.

\begin{question}\label{qu:poincare}
Does $\poi(\iot;t)$ factor through $\theta$ if and
only if $m(s,\theta(s))\in \{1,2,\infty\}$ for all $s\in S$?
\end{question}

In Section \ref{se:poincare}, in order to answer the finite case
version of Question \ref{qu:poincare} affirmatively, we resorted to a
case by case computation which did not shed much light on the
situation. A combinatorial proof, type independent if possible, would be much
preferred. Specifically, what we are looking for is this: 

\begin{problem}
Given that $\poi(\iot;t)$ factors through $\theta$, construct a bijection
$\phi:\iot \times \fix \to W$ such that $\ell(\phi((w,f))) = \rho(w) +
\ell_\theta(f)$, where $\ell_\theta$ denotes the length of
$f$ in terms of the canonical Coxeter generators of $\fix$
(cf.\ Example \ref{ex:quotient}).  
\end{problem}

Recalling from Example \ref{ex:quotient} that we may identify $\iot$
with the set of cosets $\fix \backslash W$, there is of
course a natural bijection $\iot \times \fix \to W$ defined by sending $(\fix
w,f)$ to $fw$, but in general it does not have the desired properties.

\end{document}